\documentstyle[fleqn,12pt]{article}
\begin{document}
\pagenumbering{arabic}
\setcounter{page}{1}
\pagestyle{plain}
\baselineskip=18pt

\rightline{MSUMB 95-02, April 1995} \vspace{1cm}
{\bf 
\begin{center}
QUANTUM PROPERTIES OF THE DUAL \\ MATRICES IN $GL_q(1|1)$ 
\end{center}
\vspace{2cm}
}

\begin{center} 
Salih Celik 
\end{center}
Department of Mathematics, Faculty of Sciences,
Mimar Sinan University, 80690 Besiktas, Istanbul, TURKEY.
\vspace{1cm}
\begin{center}
Sultan A. Celik 
\end{center}
Department of Mathematics, Faculty of Sciences,
Yildiz Technical University, Sisli, Istanbul, TURKEY. 

\vspace{2cm}
\begin{center}
{\bf Abstract} 
\end{center}
In this paper, we give the quantum analogue of the dual matrices for the 
quantum supergroup $GL_q(1|1)$ and discuss these properties of 
the quantum dual supermatrices.

\vfill\eject
\vspace{2cm}\noindent
{\bf 1. INTRODUCTION }

An explicit quantum deformation of the supergroup $GL(1|1)$ with two 
even and two odd generators was given by Corrigan et al in [1]. The 
properties of the 2x2-supermatrices in $GL_q(1|1)$ was investigated 
by Schwenk et al in [2]. In this work, we consider the dual supermatrices 
in $GL(1|1)$ and discuss the properties of quantum dual supermatrices.

Let us begin with some remarks. We know that the supergroup $GL(1|1)$ can be 
deformed by assuming that the linear transformations in $GL(1|1)$ are 
invariant under the action of the quantum superplane and its dual [3]. 
Consider a quantum superplane and its dual,
$$ V = \left ( \matrix{ x \cr \xi \cr } \right) ~~\mbox{and}~~ 
   \widehat{V} = \left ( \matrix{ \eta \cr y \cr} \right) \eqno(1.1)$$
satisfying 
$$ x \xi - q \xi x = 0 ,~ \xi^2 = 0, \eqno(1.2a) $$
$$ \eta^2 = 0 ,~ y \eta - q \eta y = 0 \eqno(1.2b) $$
where latin and greek letters denote even and odd elements respectively. 
Taking
$$ M = \left ( \matrix{ a & \beta \cr
                   \gamma & d \cr}
\right) \eqno(1.3) $$
as a supermatrix in $GL(1|1)$, we demand that the relations (1.2) are 
preserved under the action of $M$ on the quantum superplane $V$ and its 
dual $\widehat{V}$
$$ MV = V' ~~\mbox{and}~~ M\widehat{V} = \widehat{V}'. \eqno(1.4) $$  
We assume that even generators commute with everything and odd generators 
anticommute among themselves. Then we obtain the following $q$-commutation 
relations [1] (also see [2])
$$ a\beta = q\beta a , ~~d\beta = q\beta d, $$
$$ a\gamma = q\gamma a ,~~ d\gamma = q\gamma d, \eqno(1.5) $$
$$ \beta \gamma + \gamma \beta = 0 , ~~ \beta^2 = 0 = \gamma^2, $$
$$ ad - da = (q - q^{-1})\gamma \beta. $$
These relations will be used in sec. 4. Note that if $M \in GL_q(1|1)$ then 
$M^n \in GL_{q^n}(1|1)$. This is proved in [2].

\noindent
{\bf 2. QUANTUM DUAL SUPERMATRICES IN $GL_q(1|1)$ }

In this section we give the $q$-commutation relations which the matrix 
elements of a dual supermatrix satisfy. Let $\widehat{M}$ be a dual 
supermatrix in $GL_q(1|1)$, namely, 
$$ \widehat{M} = \left ( \matrix{ \alpha & b \cr
                                       c & \delta \cr}
\right) \eqno(2.1) $$
with its generators (anti)commuting with the coordinates of $V$ and 
$\widehat{V}$. Then, the transformations 
$$ V \longmapsto \widehat{M}V = \widehat{V}' ~~\mbox{and}~~ 
\widehat{V} \longmapsto \widehat{M}\widehat{V} = V' \eqno(2.2) $$
impose the following bilinear product relations among the generators 
of $\widehat{M}$:
  $$ \alpha b = q^{-1} b \alpha, ~~\alpha c = q^{-1}c \alpha, $$
  $$ \delta b = q^{-1}b \delta,~~ \delta c = q^{-1}c \delta, \eqno(2.3) $$
  $$ \alpha \delta + \delta \alpha = 0 , ~~ \alpha^2 = 0 = \delta^2, $$
  $$ bc - cb = (q - q^{-1})\delta \alpha $$
and $q^2 - 1 \neq 0$. 
From these relations one obtains
$$ \alpha b^{-1} = q b^{-1}\alpha, ~~ \alpha c^{-1} = q c^{-1} \alpha, $$
$$ \delta b^{-1}=qb^{-1} \delta,~~\delta c^{-1}=q c^{-1}\delta, \eqno(2.4) $$
$$ bc^{-1} - c^{-1}b = (q - q^{-1})\alpha c^{-1} \delta c^{-1} $$
provided $b$ and $c$ are invertible. 

\noindent
{\bf 3. THE INVERSE OF $\widehat{M}$}

To obtain the inverse of $\widehat{M}$, we introduce $\Delta_1$ and $\Delta_2$ 
in the form
$$ \Delta_1 = bc - q\delta \alpha ~~\mbox{and} 
~~ \Delta_2 = cb - q\alpha \delta. \eqno(3.1) $$
just as in [4]. Then one can write
$$ \widehat{M}_L^{-1} = 
  \left ( \matrix{ -q\Delta_1^{-1}\delta & \Delta_1^{-1} b \cr \cr
                     \Delta_2^{-1} c     & -q\Delta_2^{-1} \alpha \cr}
\right) \eqno(3.2) $$
as the left inverse of $\widehat{M}$. After some calculations one obtains  

$$ \Delta_1 b = b \Delta_1, ~~ \Delta_2 c = c \Delta_2, $$
$$ \Delta_k \alpha = q^2 \alpha \Delta_k, 
~~ \Delta_k \delta = q^2 \delta \Delta_k, ~~k = 1,2 \eqno(3.3) $$
and also
$$ b^2 \Delta_1^{-1} = bc^{-1} - \alpha c^{-1}\delta c^{-1}, $$
$$ c^2 \Delta_2^{-1} = cb^{-1} - \delta b^{-1}\alpha b^{-1}. \eqno(3.4)$$
Note that it is easy to verify that $b^2\Delta_1^{-1}$ and 
$c^2 \Delta_2^{-1}$ commute with everything. Therefore the matrix 
$\widehat{M}_L^{-1}$ in (3.2) may be written as 
$$ \widehat{M}_L^{-1}  
 =  \left ( \matrix{ -c^{-1}\delta c^{-1} & b^{-1} \cr \cr
                               c^{-1}     & -b^{-1}\alpha b^{-1} \cr}
\right)
\left ( \matrix{ c^2 \Delta_2^{-1} & 0 \cr \cr
                                0  & b^2 \Delta_1^{-1} \cr}
\right) \eqno(3.5) $$
which shows that $\widehat{M}_L^{-1} = \widehat{M}^{-1}$ after some  
calculations along the lines of [2], sec. 3. Thus one can define the 
quantum dual superdeterminant as follows:
$$ s\widehat{D}_q(\widehat{M}) = b^2 \Delta_1^{-1} = bc^{-1} - 
\alpha c^{-1} \delta c^{-1}. \eqno(3.6) $$
Note that the inverse of a dual supermatrix $\widehat{M}$ can be also 
obtained from the decomposition 
$$ \widehat{M} = \left(\matrix { \alpha &  b-\alpha c^{-1}\delta  \cr \cr
                                      c &     0  \cr}
\right) \left(\matrix { 1 & c^{-1} \delta  \cr \cr
                        0 & 1             \cr}
\right). \eqno(3.7) $$

Finally we note that the product of two dual supermatrices is not a dual 
supermatrix, i.e., the matrix elements of a product 
$M = \widehat{M}\widehat{M}'$ do not satisfy (2.6) but they satisfy (1.5) if 
$\widehat{M}$ and $\widehat{M}'$ are two dual supermatrices and $(b,c)$ 
($(\alpha,\delta)$) pairwise commute (anti-commute) with $(b',c')$ 
($(\alpha',\delta')$). This interesting property will show as the way to 
the contents of the next section.

\noindent
{\bf 4. PROPERTIES OF $\widehat{M}^n$}

From sec. 3 we know that the matrix elements of a product matrix 
$\widehat{M}\widehat{M}'$ obey the relations (1.5). Therefore we must 
consider the matrix elements of $\widehat{M}$ with respect to even and 
odd values of $n$. Let the $(2n-1)$-th power of $\widehat{M}$ be
$$ \widehat{M}^{2n-1} = \left(\matrix {A_{2n-1} & B_{2n-1} \cr \cr
                                       C_{2n-1} & D_{2n-1}  \cr}
\right), ~n \geq 1.  \eqno(4.1) $$
After some algebra, one obtains
$$ A_{2n-1} = \{[n]_q\alpha + q[n-1]_q\delta\}(bc)^{n-1}, $$
$$ B_{2n-1} = \{bc + q[n-1]_{q^2}\alpha\delta\}(bc)^{n-2}b,  \eqno(4.2) $$
$$ C_{2n-1} = \{cb + q[n-1]_{q^2}\delta\alpha\}(cb)^{n-2}c,  $$
$$ D_{2n-1} = \{[n]_q\delta + q[n-1]_q\alpha\}(cb)^{n-1},  $$
where
$$ [n]_q = \frac{1 - q^{2n}}{1-q^2}. \eqno(4.3) $$
Now it is easy to show that the following relations are satisfied.
$$ A_{2n-1} B_{2n-1} = q^{-(2n-1)}B_{2n-1}A_{2n-1} $$
$$ A_{2n-1} C_{2n-1} = q^{-(2n-1)}C_{2n-1}A_{2n-1} $$
$$ D_{2n-1} B_{2n-1} = q^{-(2n-1)}B_{2n-1}D_{2n-1} $$
$$ D_{2n-1} C_{2n-1} = q^{-(2n-1)}C_{2n-1}D_{2n-1}, \eqno(4.4) $$
$$ A_{2n-1} D_{2n-1} + D_{2n-1} A_{2n-1} = 0,$$
$$A_{2n-1}^2 = 0 = D_{2n-1}^2,  $$
$$ B_{2n-1} C_{2n-1} - C_{2n-1}B_{2n-1} = (q^{2n-1} - q^{-(2n-1)})
A_{2n-1}D_{2n-1}. $$
Then $\widehat{M}^{2n-1}$ is a dual supermatrix with deformation parameter 
$q^{2n-1}$.

Similarly, if we write for the matrix $\widehat{M}^{2n}$, the (2n)-th 
power of $\widehat{M}$ as 
$$ \widehat{M}^{2n} = \left(\matrix {A_{2n} & B_{2n} \cr \cr
                                       C_{2n} & D_{2n}  \cr}
\right), ~n \geq 1  \eqno(4.5) $$
where (after some calculations)

$$ A_{2n} = \{bc + q\frac{1-q^2}{1+q^2}[n]_q[n-1]_q\alpha\delta\}(bc)^{n-1},$$
$$ B_{2n} = [n]_q\{\alpha + q\delta\}b(cb)^{n-1},  \eqno(4.6) $$
$$ C_{2n} = [n]_q\{\delta + q\alpha\}c(bc)^{n-1},  $$
$$ D_{2n} = \{bc + q\frac{1-q^2}{1+q^2}[n]_q[n-1]_q\delta\alpha\}(cb)^{n-1}, $$
then the elements of $\widehat{M}^{2n}$ obey the following relations
$$ A_{2n} B_{2n} = q^{2n}B_{2n}A_{2n} $$
$$ A_{2n} C_{2n} = q^{2n}C_{2n}A_{2n} $$
$$ D_{2n} B_{2n} = q^{2n}B_{2n}D_{2n} $$
$$ D_{2n} C_{2n} = q^{2n}C_{2n}D_{2n},\eqno(4.7) $$ 
$$ B_{2n} C_{2n} + C_{2n} B_{2n} =0, $$
$$B_{2n}^2 = 0 = C_{2n}^2,  $$
$$ A_{2n} D_{2n} - D_{2n}A_{2n} = (q^{2n} - q^{-2n})
C_{2n}B_{2n}. $$
Thus the matrix $\widehat{M}^{2n}$ is a supermatrix in the form (1.3).

Equations (4.4) and (4.7) can be proved using the relation (2.3).

\noindent
{\bf 5. CONCLUSIONS} 

We have given the $q$-commutation relations which the matrix elements of a 
dual supermatrix in $GL_q(1|1)$ satisfy and obtained the (dual) quantum 
super-inverse and (dual) quantum superdeterminant of a dual quantum 
supermatrix. Finally we have shown that it must consider the matrix 
elements of a dual supermatrix with respect to even and odd values of $n$. 
And so we discussed the properties of the $n$-th power of a dual supermatrix.

{\bf Acknowledgment.} We wish to thank Prof. Dr. M. Arik for reading the 
manuscript and for many useful suggestions and clarifications.

{\bf REFERENCES }

\begin{description}
\item[[1]] Corrigan, E., Fairlie, B., Fletcher, P. and Sasaki, R., 
           J. Math. Phys. {\bf 31}, 776, 1990.
\item[[2]] Schwenk, J., Schmidke, B. and Vokos, S., 
           Z. Phys. C {\bf 46}, 643, 1990.
\item[[3]] Manin, Yu I., 
           Commun. Math. Phys. {\bf 123}, 163, 1989.
\item[[4]] Celik, S., Celik, S. A., 
           {\it On the quantum supergroup $SU_{p,q}(1|1)$ and quantum 
            oscillators}, Preprint MSUMB - 95/1, 1995. 
\end{description}

\end{document}